\documentclass[12pt]{article}
\usepackage{amsfonts,amssymb,version}

\def\AA{{\bf A}}

\def\GG{{\mathbb G}}

\def\HNP{\mathop{\rm HNP}\nolimits}
\def\hn{\mathop{\rm HN}\nolimits}
\def\HNT{\mathop{\rm{HNT}}\nolimits}
\def\hnt{\mathop{\rm{HNT}}\nolimits}
\def\HOM{\mathop{\rm Hom}\nolimits}

\def\OO{{\mathcal O}}

\def\Q{\mathbb Q}

\def\Z{\mathbb Z}

\def\deg{\mathop{\rm deg}\nolimits} 
\def\gr{\mathop{\rm Gr}\nolimits}

\def\ker{\mathop{\rm ker}\nolimits}
\def\qed{\hfill$\square$\medskip}

\def\spec{\mathop{\rm Spec}\nolimits}

\def\toto{\stackrel{\to}{{\scriptstyle \to}}}

\let\hra\hookrightarrow
\let\ov\overline
\let\wh\widehat

\newtheorem{theorem}{Theorem}[section]
\newtheorem{proposition}[theorem]{Proposition}
\newtheorem{lemma}[theorem]{Lemma}

\newtheorem{corollary}[theorem]{Corollary}

\def\rem{\refstepcounter{theorem}\paragraph{Remark \thetheorem}}

\def\proof{\paragraph{Proof}}

\makeatletter \def\l@section{\@dottedtocline{1}{0em}{1.2em}} \makeatother

\begin{document}

\centerline{\Large\bf Schematic HN stratification for families}

\centerline{\Large\bf of principal bundles and lambda modules}

\bigskip

\centerline{\bf Sudarshan Gurjar and Nitin Nitsure} 


\begin{abstract}
For a family of principal bundles with a reductive structure group
on a family of curves in characteristic zero, it is known that the 
Harder Narasimhan type of its restriction to each fiber varies
semicontinuously over the parameter scheme of the family. This 
defines a stratification of the parameter scheme by locally closed
subsets, known as the Harder-Narasimhan stratification.

In this note, we show how to endow each Harder-Narasimhan stratum 
with the structure of a locally closed subscheme of the parameter
scheme, which enjoys the universal property that under any base change 
the pullback family admits a relative Harder-Narasimhan filtration
with a given Harder-Narasimhan type if and only if the base change 
factors through the  schematic stratum corresponding to that 
Harder-Narasimhan type. This has the consequence that principal
bundles of a given Harder Narasimhan type form an Artin stack.

We also prove a similar result showing the
existence of a schematic Harder-Narasimhan filtration for
flat families of pure sheaves of $\Lambda$-modules (in the sense of Simpson) 
in arbitrary dimensions
and in mixed characteristic, generalizing the result for 
sheaves of ${\mathcal O}$-modules proved earlier by Nitsure. 
This again has the 
implication that $\Lambda$-modules of a fixed 
Harder-Narasimhan type form an Artin stack. 
\end{abstract}


\centerline{2010 Math. Subj. Class. : 14D20, 14D23, 14F10.}

\section*{Introduction}

Let $G$ be a reductive algebraic group over a field $k$ of characteristic zero,
and let $T\subset B \subset G$ be a chosen maximal torus and a 
Borel subgroup of $G$. We assume that $T$ is split over $k$. 
If $X$ is a curve over $k$ (always assumed to be smooth and projective over
$k$ and geometrically irreducible)
and if $E$ is a principal $G$-bundle on $X$ 
(locally trivial in \'etale topology), 
then it is known that $E$ admits a unique 
{\it canonical reduction} $(P,\sigma)$
(see Atiyah-Bott [A-B], Behrend [Be] and Biswas-Holla [Bi-Ho]). Here, 
$B\subset P\subset G$ is a standard parabolic subgroup of $G$ 
and $\sigma : X\to E/P$ is a section which gives a reduction 
of structure group of $E$ from $G$ to $P$. Any reduction 
$(P',\sigma')$ of $E$ to an arbitrary parabolic subgroup 
$B\subset P'\subset G$ gives rise to a point of the closed 
positive Weyl chamber 
$\ov{C} \subset \Q\otimes X_*(T)$. The particular such point  
defined by the canonical reduction of $E$ is called the
{\it Harder Narasimhan type} (or simply the HN type) 
of $E$. We denoted it by $\hn(E) \in \ov{C}$.
Given any family of curves
$X/S$, where $S$ is a noetherian $k$-scheme, 
and a family of principal bundles $E$ over it, the 
HN type $\hn(E_s)$ of its restriction $E_s = E|_{X_s}$ to the fiber 
$X_s$ of $X$ over $s\in S$ is known to be an upper semicontinuous function on 
the underlying topological space $|S|$ of the parameter scheme $S$
(we give a proof later, as it is not easy to find a complete proof 
in the literature)
with respect to a natural partial order on the closed positive Weyl chamber.
Thus, for each point $\tau \in \ov{C}$, the level 
subset $|S|^{\tau}(E)$ of $|S|$ where $\hn(E_s) = \tau$
is locally closed. 

Our main result is that for any family $E$ over $X/S$, each such 
level set $|S|^{\tau}(E)$ admits a unique structure of a locally closed 
subscheme $S^{\tau}(E)\subset S$, which has the following universal property.
If $T\to S$ is a morphism of $k$-schemes, such that the base-change
$E_T$ which is a principal bundle over $X_T$ admits a {\it relative
canonical reduction} over $T$ with constant HN type $\tau$, then 
$T\to S$ factors uniquely via $S^{\tau}(E)\hra S$.

This note is organized as follows. In section 1, we quickly recall all
basic definitions and results about principal bundles that we need.

In section 2, we prove our main result and some corollaries. 

Section 3 is devoted to the HN stratification of the parameter 
scheme $S$ for a family $E$ of $\OO$-coherent, flat, pure dimensional
$\Lambda$-modules (in the sense of Simpson [Si]) 
over a family $X/S$ of projective schemes over an 
arbitrary locally noetherian scheme $S$ (which may have mixed characteristic). 
Generalizing the main result of Nitsure [Ni-3], which is the case when 
$\Lambda = \OO$, we show that 
there exists a schematic Harder Narasimhan stratification 
$\cup_{\tau}\, S^{\tau}(E)$ of $S$, with 
the universal property that an arbitrary base change $T\to S$ 
factors via a stratum $S^{\tau}(E)$ if and only if the base-change 
$E_T$ admits a relative HN filtration (as $\Lambda$-modules) of 
constant HN type $\tau$ on $T$.


\pagestyle{myheadings}
\markright{Gurjar and Nitsure: HN stratifications for principal bundles
and $\Lambda$-modules}


\section{Preliminaries on principal bundles}

Let $k$ be a field of characteristic zero, and let 
$G\supset B \supset T$ be a 
reductive group, together with a chosen Borel subgroup and a maximal 
torus (we assume that $G$ is split over $k$). 
As usual, $X^*(T)$ and $X_*(T)$ will respectively denote the groups of all 
characters and all $1$-parameter subgroups of $T$.
Let $\Delta \subset X^*(T)$ be the corresponding simple roots.
Let $\omega_{\alpha} \in \Q\otimes X^*(T)$ denote the
fundamental dominant weight corresponding to $\alpha \in \Delta$,
so that $\langle \omega_{\alpha}, \beta^{\vee}\rangle = \delta_{\alpha, \beta}$
where $\beta^{\vee} \in \Q\otimes X_*(T)$ is the simple coroot corresponding 
to $\beta \in \Delta$. Note that each $\omega_{\alpha}$ is a non-negative
linear combination (with coefficients in $\Q^{\ge 0}$) 
of the simple roots $\alpha$.
Recall that {\it the closed positive Weyl chamber
$\ov{C}$} is the subset of $\Q\otimes X_*(T)$
defined by the condition that
$\mu \in \ov{C}$ if and only if  
$\langle \alpha, \mu\rangle \ge 0$ for all $\alpha \in \Delta$.  
The standard {\it partial order} on $\ov{C}$  
is defined by putting $\mu \le \pi$ 
if $\omega_{\alpha}(\mu) \le \omega_{\alpha}(\pi)$ for all $\alpha \in \Delta$,
and $\chi(\mu) = \chi(\pi)$ for all $\chi: G \to \GG_m$.

By definition, all roots and weights in $X^*(T)$ are trivial on
the connected component $Z_0(G)\subset T$ of the center of $G$. 

We choose, once for all, a Weyl group invariant positive definite 
bilinear form on $\Q\otimes X^*(T)$ taking values in $\Q$. This, 
in particular, will allow us to identify $\Q\otimes X^*(T)$ with 
$\Q\otimes X_*(T)$.

Let $I\subset \Delta$ be any subset.
The following elementary property (which occurs in 
[Bi-Ho], and which holds in any abstract root system simply
because $\oplus_{\beta\in I}\Q\beta$ is the orthogonal complement of 
$\oplus_{\gamma \in \Delta -I}\Q\omega_{\gamma}$ and 
the angle between any two simple roots is $\ge \pi/2$ radians) 
will be useful in what follows. 

(R) Given any $\alpha \in \Delta -I$, there exist 
(unique) rational numbers $n_{\beta}\ge 0$ where $\beta\in I$, 
such that the element 
$\chi_{\alpha} = \alpha + \sum_{\beta \in I} n_{\beta}\beta$
lies in $\oplus_{\gamma \in \Delta -I} \Q\omega_{\gamma}$.

\bigskip

\centerline{\bf Canonical reductions of principal bundles}

\medskip

Let $T\subset B\subset G$ be as above. Recall that 
a principal $G$-bundle $E$ on a smooth irreducible projective 
curve $X/k$ is said to be {\it semistable} if for any 
reduction $\sigma : X\to E/P$ of the structure group to a parabolic $P\subset G$
and any dominant character $\chi: P\to \GG_m$, we have 
$$\deg(\chi_*\sigma^*E) \le 0$$
where $\sigma^*E$ is the principal $P$-bundle on $X$ 
defined by the reduction $\sigma$,
and $\chi_*\sigma^*E$ is the $\GG_m$-bundle obtained by
extending its structure group via $\chi: P \to \GG_m$,
which is equivalent to a line bundle on $X$.

A {\it canonical reduction} of a principal $G$-bundle $E$ 
is a pair $(P,\sigma)$ where $P$ is a standard parabolic subgroup
of $G$ (that is, $B\subset P$ for the chosen Borel $B$) and 
$\sigma : X\to E/P$ is a reduction of the structure group to $P$ 
such that the following two conditions (C-1) and (C-2) hold. 

{\bf (C-1)} If $\rho: P \to L = P/U$ is the Levi quotient of $P$ (where $U$ is
the unipotent radical of $P$) then the principal $L$-bundle 
$\rho_*\sigma^*E$ is semistable.

{\bf (C-2)} For any non-trivial character $\chi: P \to \GG_m$ 
whose restriction to the chosen maximal torus $T\subset B \subset P$ 
is a non-negative linear combination $\sum n_i\alpha_i$ 
of simple roots $\alpha_i \in \Delta$ (where $n_i \ge 0$, and at least
one $n_i \ne 0$), we have $\deg(\chi_*\sigma^*E) > 0$.

\rem\label{can red}
It was shown by Behrend [Be] (see Biswas-Holla [Bi-Ho] 
for another proof) that each $E$ has a unique canonical reduction. 
Moreover, under any extension of base fields, it is easy to see
that the canonical reduction base changes.

The standard parabolics $G\supset P \supset B$ of $G$ are in one-one 
correspondence with the subsets of $\Delta$, under which to any $P$
there is associated the subset $I_P\subset \Delta$ of the corresponding
inverted simple roots (for example, $I_G = \Delta$ while $I_B = \emptyset$).
Let $\wh{P}$ denote the group of all characters $\chi : P \to \GG_m$,
and let $\wh{P}|_T\subset X^*(T)$ be the group of 
their restrictions to $T$. 
We have an internal direct sum decomposition
$$\Q\otimes X^*(T) = (\Q\otimes \wh{P}|_T) \oplus 
(\oplus_{\alpha \in I_P} \Q \alpha),$$
where, in turn, 
$$\Q\otimes \wh{P}|_T = (\oplus_{\alpha \in \Delta - I_P}
\Q\omega_{\alpha}) \oplus (\Q\otimes \wh{G}|_T).$$

Given a reduction $(P,\sigma)$ of $E$, we get 
an element $\mu_{(P,\sigma)}\in \Q\otimes X_*(T)$ defined by 
\begin{eqnarray*}
\langle \chi, \mu_{(P,\sigma)}\rangle & = & \left\{
\begin{array}{ll}
\deg(\chi_*\sigma^*E) & \mbox{if }  \chi \in \wh{P},  \\
0                     & \mbox{if } \chi \in I_P.
\end{array}\right.
\end{eqnarray*}

If $(P,\sigma)$ is the canonical reduction of $E$, then the element
$\mu_{(P,\sigma)}$ is called the {\it HN type of $E$}, and is denoted by 
$\hn(E)$. If $\alpha \in \Delta -I_P$ and 
$\chi_{\alpha} 
= \alpha + \sum_{\beta \in I_P} n_{\beta}\beta \in 
\oplus_{\gamma \in \Delta -I_P} \Q\omega_{\gamma}$, 
where $n_{\beta}\in \Q^{\ge 0}$, 
is the rational character given by the property (R) of root systems 
mentioned earlier, then 
$$\langle \chi_{\alpha}, \hn(E)\rangle = 
\deg({\chi_{\alpha}}_*\sigma^*E) \ge 0$$
for all $\alpha \in \Delta -I_P$ by the condition (C-2) above.
As $\langle \beta,\hn(E)\rangle = 0$ for all $\beta \in I_P$,
we get $\langle \alpha,\hn(E)\rangle = \langle \chi_{\alpha}, \hn(E)\rangle 
\ge 0$ for all $\alpha \in \Delta -I_P$. Hence 
$\hn(E)$ is in the closed positive Weyl chamber $\ov{C}$, in fact, 
in the facet of $\ov{C}$ defined by the vanishing of all $\beta \in I_P$.

Note that a principal bundle $E$ of type $\hn(E)= \mu$ is semistable
if and only if $\mu$ is central, that is, $\mu = a\nu$ for
some $1$-parameter subgroup $\nu : \GG_m \to Z_0(G)$ and $a\in \Q$.  

Given the HN-type $\mu = \hn(E)$ of $E$, we can recover the corresponding
standard parabolic $P$ as follows. Let $I_{\mu} \subset \Delta$ 
be the set of all simple roots $\beta$ such that $\langle \beta ,\mu\rangle =0$.
Then $I_{\mu}$ is exactly the set of inverted simple roots which defines $P$. 
Alternatively, let $n \ge 1$ be any integer such that $\nu = n \mu \in X_*(T)$.
The the $k$-valued points of $P$ are all those $g$ for which 
$\lim_{t\to 0} \nu(t)g\nu(t)^{-1}$ exists in $G$.

For any standard parabolic $P$, let 
$\theta_P = -\sum_{\Delta- I_P}\alpha : P\to \GG_m$ be the character
which defines the ample line bundle $\det(T_{G/P})$ on $G/P$. 
As remarked in [Bi-Ho], for any reduction $(P,\sigma)$ of a given 
$E$, the vector bundle $\sigma^*(T_{E/P})$ is a quotient of 
$Ad(E)$ and has determinant equal to 
${\theta_P}_*\sigma^*E$, hence the set of 
integers $\deg {\theta_P}_*\sigma^*E$ as $(P,\sigma)$
varies is bounded below. Let its minimum element be denoted by $d_E$. 

We will find the following fact useful (see [Bi-Ho] Remark 3.2):

{\bf (C-3)} A parabolic reduction 
$(P,\sigma)$ is the canonical reduction of $E$ if and only if 
$\deg {\theta_P}_*\sigma^*E = d_E$
and (C-2) is satisfied. 

Therefore by the definition of $\hn(E)$, 
when $(P,\sigma)$ is the canonical reduction of $E$, we get 
$\langle \theta_P , \hn(E)\rangle =d_E$.

We will also use the following elementary, well-known fact.

{\bf (C-4)} Let $E$ be a principal $G$-bundle on a curve $X$, let
$(P,\sigma)$ be its canonical reduction, and let $(Q,\tau)$ be any
reduction to a standard parabolic. Let $\hn(E) = \mu_{(P,\sigma)}$
and $\mu_{(Q,\tau)}$ be the corresponding elements of $\Q\otimes X_*(T)$
(the element $\hn(E)$ lies in the closed positive Weyl chamber $\ov{C}$, 
but $\mu_{(Q,\tau)}$ need not do so). Then for each $\alpha \in \Delta$ 
we have the inequality
$$\langle \omega_{\alpha}, \hn(E)\rangle = 
\langle \omega_{\alpha}, \mu_{(P,\sigma)}\rangle \ge  
\langle \omega_{\alpha}, \mu_{(Q,\tau)}\rangle $$
where $\omega_{\alpha} \in \Q\otimes X^*(T)$ is the fundamental dominant weight
corresponding to $\alpha$.
Moreover, if each of the above inequalities is an equality then 
$(P,\sigma) = (Q,\tau)$.

\bigskip

\centerline{\bf Families of principal bundles}

We now fix a smooth projective morphism $X\to S$ of relative dimension one, 
such that each geometric fiber is irreducible, and
$S$ is a noetherian scheme over $k$. By a {\it family} of 
principal $G$-bundles on $X/S$ parameterized by an $S$-scheme $T$
one means a principal $G$-bundle $E$ on $X_T = X\times_ST$. The
following facts are well-known for any such family.

{\bf (1) Openness of semistability.} All $t\in T$ such that $E_t = E|_{X_t}$ is
semistable form an open subset of $|T|$ (see Ramanathan [R] Corollary 3.18).
The corresponding open subscheme of $T$ will be called the open subscheme
of semistable bundles.

{\bf (2) Semicontinuity of HN type.} Given any $\mu \in \ov{C}$,
the set $|T|^{\mu}$ 
of all $t\in T$ such that $E_t = E|_{X_t}$ is of type $\mu$ is a locally
closed subset of $|T|$. The closure of any $|T|^{\mu}$ is contained 
in the union of all $|T|^{\nu}$ for $\nu \ge \mu$. 
This is well-known (a sketch of a proof that (1) implies (2) is given in the 
next section).

Let $E$ be a family of principal $G$-bundles on $X/S$, 
such that the HN-type $\hn(E_s)$ is constant (say $= \tau \in \ov{C}$)
for each $s\in S$. 
A {\it relative canonical reduction} for $E/X/S$
is a reduction $(P,\sigma)$ of the structure group of $E$ to 
a standard parabolic $P$ over all of $X$, 
which restricts to give the canonical reduction of $E_s$ over $X_s$
for each $s\in S$.
Note that even if $\hn(E_s)$ is constant, 
in general such a relative canonical reduction need not exist.
However, as we will prove (see Theorem \ref{six} and 
Corollary \ref{2.7} below), it is unique whenever it exists, and it 
exists in particular when $S$ is reduced.

\section{HN stratification for families of principal bundles}

We begin by recalling two results of Grothendieck [FGA] on Hilbert schemes 
(see for example [Ni-1] for an exposition).

(A) Given a projective morphism  $W \to S$ 
there exists a relative Hilbert scheme 
$Hilb_{W/S} \to S$ which represents the contravariant 
functor from $S$-schemes to $Sets$ which to any 
$S$-scheme $T\to S$ associates the set of all closed 
subschemes $Y\subset W_T = W\times_ST$ such that $Y\to T$ is 
flat. The scheme $Hilb_{W/S}$ is a disjoint union of open (and closed) 
subschemes which are projective over $S$, in particular,
the morphism $Hilb_{W/S} \to S$ satisfies the valuative criterion for
properness.  

(B) Given any projective morphisms $W \to X \to S$ where
$X\to S$ is flat, 
there is an open subscheme $R_{W/X/S}\subset Hilb_{W/S}$ 
which parameterizes all sections of $W \to X$ on 
fibers of $X\to S$. It represents the contravariant 
functor from $S$-schemes to $Sets$ which associates
to any $S$-scheme $T\to S$ the set of all sections of $W_T \to X_T$. 

For any parabolic subgroup $P\subset G$, and a family $E$ of principal 
$G$-bundles on $X/S$, 
the projection $\pi : E/P \to X$ is a smooth projective morphism. 
Any character $\chi : P\to \GG_m$ defines a line bundle $L_{\chi}$
on $E/P$ which is relatively ample over $X$ whenever $\chi$ is 
negative dominant. By the result (A) above, 
we get a relative Hilbert scheme 
$$Hilb_{(E/P)/S} \to S$$ 
which parameterizes all closed subschemes of the fibers of 
$E/P \to S$. By the result (B) above, this has an open subscheme 
$$R_{(E/P)/X/S} \subset Hilb_{(E/P)/S}$$
which parameterizes all sections of $E/P \to X$ on 
fibers of $X\to S$. It represents the contravariant 
functor from $S$-schemes to $Sets$ which associates
to any $S$-scheme $T\to S$ the set of all reductions 
$(P, \sigma: X_T \to E_T/P)$
of structure group of $E_T$ from $G$ to $P$.

\begin{lemma}\label{one} Let $X/S$ be a smooth projective family of 
geometrically irreducible curves over a noetherian integral scheme $S$ over $k$.
Let $K$ be the function field of $S$. Let  
$E$ be a principal $G$-bundle on $X$, let $P\subset G$ be a standard 
parabolic, and let 
$(P, \tau : X_K \to E_K/P)$ be a canonical reduction of $E_K$. 
Then there exists a non-empty open subscheme $U\subset S$ 
and a section $\sigma: X_U \to E_U/P$ such that 
$(P,\sigma)$ is a relative canonical reduction of $E_U/X_U/U$.
\end{lemma}

\proof The reduction $(P, \tau)$ defines a $K$-valued point 
$\tau \in  R_{(E/P)/X/S} $. Under the inclusion 
$R_{(E/P)/X/S} \subset Hilb_{(E/P)/S}$, let this point
lie in an open and closed subscheme $H^o \subset Hilb_{(E/P)/S}$
which is projective over $S$. Let $R^o = R_{(E/P)/X/S}\cap H^o$,
so $\tau \in  R^o$. As $H^o\to S$ is projective, there exists 
a nonempty open subscheme $U_1 \subset S$ and a section $\tau_1 : U_1 \to H^o$
which extends $\tau$. Let $U_2 = \tau_1^{-1}(R^o)\subset S$, 
which is open and is nonempty as it contains the generic point $\spec K$. 
This gives a reduction $\tau_2 :  X_{U_2} \to E_{U_2}/P$ which extends
$\tau : X_K \to E_K/P$. Its associated Levi bundle is semistable
over the generic point of $U_2$, so by openness of semistability,
the associated Levi bundle 
is semistable over an open neighbourhood $U$ of the generic point in $U_2$.
Taking $\sigma = \tau_2|_U$, we get the desired 
relative canonical reduction $(P,\sigma)$ of $E_U$ over $X_U$.
The condition (C-1) is satisfied by the definition of $U$, and 
the condition(C-2) is satisfied by the constancy of 
$\deg \chi_*\sigma_s^*E_s$ as $s$ varies over $U$. \qed

\begin{lemma}\label{two} Let $X$ be a noetherian integral scheme of 
dimension $\ge 1$ which is normal.
Let $\pi: W \to X$ be a projective morphism.
Let $F\subset X$ be  a closed subset of dimension zero and  
let $\sigma : X-F \to W$ be a section  of $\pi$ outside $F$.
Let $V\subset W$ be the closure of $\sigma(X-F)$. If
for each $x\in F$, the set $\pi^{-1}(x)\cap V$ is finite
then $\sigma$ prolongs to a global section of $\pi: W \to X$.
\end{lemma}

\proof It is enough to assume that $F$ is a singleton set $\{ x\}$
for some closed point $x\in X$.
We can take an affine open $U$ in $W$ such that the 
finite set $\pi^{-1}(x)\cap V$ is inside $U$. Then by properness of 
$\pi$, the subset $X' = X - \pi(V-U)$ is open in $X$, and $F\subset X'$. 
The restricted section $\sigma|_{X' -F} : X' -F \to U$ is given by
regular functions as $U$ is affine, so prolongs to 
$F$ by the normality of $X'$.  \qed

\begin{lemma}\label{three} Let $R$ be a discrete valuation ring over $k$, 
with quotient field $K$ and residue field $k_1$,  
let $X$ be a smooth projective family of curves over $S = \spec R$ 
with geometrically 
irreducible fibers $X_0 = X\otimes_RK$ and $X_1= X\otimes_Rk_1$
over the generic point $\eta_0 = \spec K$ and the special point 
$\eta_1 = \spec k_1$ of $S$. Let $E$ be a principal $G$-bundle on 
$X$. Then the following holds.

(1) Between the HN-types of the restrictions $E_0 = E|_{X_0}$ and 
$E_1 = E|_{X_1}$, we have the inequality $\hn(E_0) \le \hn(E_1)$.

(2) If $\hn(E_0) = \hn(E_1)$, then there exists a relative 
canonical reduction for $E/X/S$.
\end{lemma}

\proof (1) Let $\alpha : X_0 \to E_0/P_0$ be the canonical reductions
of $E_0$, where $P_0$ is a standard parabolic in $G$. 
Let $Hilb_{(E/P_0)/R}$ denote the relative 
Hilbert scheme which parameterizes all closed subschemes of
the fibers of $E/P_0 \to \spec R$. The image $\alpha(X_0)$ 
of $\alpha$ defines a valued point
$$[\alpha(X_0)] \in Hilb_{(E/P_0)/R}(K).$$ 
By the valuative criterion of properness, which is satisfied by 
$Hilb_{(E/P_0)/R} \to \spec R$, the valued point
$[\alpha(X_0)] \in Hilb_{(E/P_0)/R}(K)$ prolongs uniquely to a valued point
$$ [Y]  \in Hilb_{(E/P_0)/R}(R)$$
which is represented by a closed subscheme $Y\subset E/P_0$ which is flat 
over $R$, and whose fiber $Y_0$ over $\eta_0 = \spec K$ is $\alpha(X_0)$.

The fiber $Y_1$ of $Y$ over the special point $\eta_1 = \spec (k_1)$
of $S$ is of dimension $1$. The projection $\pi_1: Y_1\to X_1$
is surjective as the projection $\pi : Y \to X$ is generically
surjective over $X$, for $\pi_0 : \alpha(X_0) \to X_0$ is surjective,
and $Y_0 = \alpha(X_0)$. Note that $X$ is $2$-dimensional and regular.
Hence by the valuative criterion of 
properness applied to the projection $\pi: E/P_0 \to X$, 
the section $\alpha: X_0 \to E_0/P_0$ prolongs to a section
$$\sigma : X - F \to E/P_0 - \pi^{-1}(F)$$
where $F \subset X_1$ is a finite set of closed points of $X_1$.
By closedness of $Y$, we must have
$\sigma(X - F) \subset Y$. 
By the valuative criterion of 
properness applied to the projection $\pi_1: E_1/P_0 \to X_1$, 
the section 
$\sigma_1 = \sigma|_{X_1-F} : X_1 -F  \to E_1/P_0$
prolongs to a section
$\gamma : X_1 \to E_1/P_0$.
Then 
$\gamma(X_1)\subset Y_1$. 
But $Y_1$ is of dimension $1$, so $\gamma(X_1)$ is an 
irreducible component of $Y_1$. Let $Z$ be the union of the remaining
irreducible components of $Y_1$, if any. Hence we get
$$Y_1 = \gamma(X_1) \cup Z.$$

If $\delta_i : X_i \to E_i/P_i$ is any reduction, then the 
Hilbert polynomial of $[\delta_i(X_i)]$ w.r.t. any line bundle $L$ on $E/P$ is 
the Hilbert polynomial of $X_i$ w.r.t. $\delta_i^*L$, that is,
$$p_L(\delta_i(X_i))(n) = \dim H^0(X_i, \delta_i^*L^{\otimes n}) 
\mbox{ for } n >> 0.$$
If $g$ denotes the genus of $X_i$, then by Riemann-Roch this gives 
$$p_L(\delta_i(X_i))(n) = (\deg \delta_i^*L)n + 1 -g.$$

By flatness of $Y$ over $R$, the Hilbert polynomial $p_L([Y_1])$ 
of $[Y_1]$ is the Hilbert polynomial $p_L([\alpha(X_0)])$ of 
$[\alpha(X_0)]$. As a non-empty $Z$ will a 
positive Hilbert polynomial whenever $L$ is relatively ample on 
$E_1/P_0 \to X_1$, from $Y_1 = \gamma(X_1) \cup Z$ 
we now get the following crucial fact. 

{\bf (*)} If a line bundle $L_{\chi}$ on $E/P_0$ is defined by
a negatively dominant character $\chi : P_0 \to \GG_m$ 
(so it is relatively ample over $X$) then 
we have the inequality of Hilbert polynomials 
$p_{L_{\chi}}([\gamma(X_1)]) \le p_{L_{\chi}}([\alpha(X_0)])$ 
which is an equality if and only if $Z$ is empty. 
Equivalently by Riemann-Roch, we get
$$\deg \chi_*\gamma^*E_1 \le \deg \chi_*\alpha^*E_0$$
where equality holds if and only if $Z = \emptyset$.

By the definition of the element $\mu_{(P_0,\sigma)} \in \Q\otimes X_*(T)$ 
introduced in the course of defining the HN-type of a principal bundle, 
this means
$$\langle \chi, \mu_{(P_0,\gamma)}\rangle \le 
 \langle \chi, \hn(E_0) \rangle \mbox{ for all negative dominant }
\chi : P_0 \to \GG_m.$$

Replacing negative dominant characters by fundamental positive 
dominant characters $\omega_{\alpha}$ 
(which changes the direction of the inequality), we get the following. 

{\bf (**)} If $I_0$ is the set of inverted simple roots defining $P_0$, 
we have
$$\langle \omega_{\alpha}, \mu_{(P_0,\gamma)}\rangle \ge 
 \langle \omega_{\alpha}, \hn(E_0) \rangle \mbox{ for all }
\alpha \in \Delta - I_0.$$

By property (C-4) of canonical reductions we have 
$$\langle \omega_{\alpha}, \hn(E_1) \rangle \ge 
\langle \omega_{\alpha}, \mu_{(P_0,\gamma)}\rangle 
\mbox{ for all }\alpha \in \Delta.$$

Hence we get
$$\langle \omega_{\alpha}, \hn(E_1)  \rangle \ge  
\langle \omega_{\alpha}, \hn(E_0) \rangle 
\mbox{ for all } \alpha \in \Delta - I_0.$$
Note that $\hn(E_0)$ lies in the facet $C(I_0) \subset \ov{C}$ 
defined by the vanishing of all $\beta \in I_0$.
Hence the above inequalities imply that $\hn(E_1) \ge \hn(E_0)$,
which proves (1).

Next we prove (2). As by assumption $\hn(E_0) = \hn(E_1)$, 
there is a common standard parabolic $P = P_0\subset G$, and canonical
reductions $(P,\alpha)$ for $E_0$ and $(P, \beta)$ for $E_1$.
We must show that there exists a section $\sigma : X \to E/P$
such that $\sigma_0 =  \alpha$ and  $\sigma_1 =  \beta$
where $\sigma_i$ denotes $\sigma|_{X_i}$.

Consider the relatively ample line bundle $L$ on $E/P$ over $X$
which comes from the negative dominant character $\theta_P$ on $P$. 
By the property (C-3) of a canonical reduction we must have 
$d_{E_1} \le \deg {\theta_P}_*\gamma^*E_1$ and 
$\deg {\theta_P}_*\alpha^*E_0 = d_{E_0}$. Hence from (*) above 
we get 
$$d_{E_1} \le \deg {\theta_P}_*\gamma^*E_1 \le 
\deg {\theta_P}_*\alpha^*E_0  = d_{E_0}.$$
But as by assumption $\hn(E_0) = \hn(E_1)$, we have $d_{E_0}= d_{E_1}$,
so we get the equality 
$$\deg {\theta_P}_*\gamma^*E_1 = 
\deg {\theta_P}_*\alpha^*E_0.$$

Hence we can conclude by (*) that $Z$ is empty, so $Y_1 = \gamma(X_1)$.
As $Y_1 = \gamma(X_1)$, the fibers of $Y$ over $F$ are singletons. 
By Lemma \ref{two}, this implies that $\sigma : X-F \to E/P$ 
prolongs to a global section over $X$, 
so we can take $F$ to be empty, and regard $\sigma : X\to E/P$
as a global section. 
As $\gamma$ is the restriction $\sigma|_{X_1}$, 
the reduction $(P,\gamma)$ of $E_1$ satisfies (C-2) by continuity from
$\alpha = \sigma|_{X_0}$. 
As $d_{E_1} = \deg {\theta_P}_*\gamma^*E_1$, it follows by (C-3) that
$(P,\gamma)$ is the canonical reduction of $E_1$. 
Thus, $\sigma : X  \to E/P$ 
is the desired  relative canonical reduction
for $E$ over $R$.  This completes the proof of (2), and of Lemma \ref{three}.
\qed

The lemmas \ref{one} and \ref{three}.(1) 
allow us to give a proof of the following 
known result.

\begin{proposition}\label{four} Let $X$ be a smooth family of 
geometrically irreducible curves over locally noetherian scheme $S$ over $k$.
Let $E$ be a principal $G$-bundle on $X$. Then the function
$|S| \to \ov{C}$ under which $s\mapsto \hn(E_s)$ is upper
semicontinuous.  Consequently, for each $\tau \in \ov{C}$, 
the level subset $|S|^{\tau}(E)$ is locally closed in $|S|$, 
and for any maximal element $\tau \in \{ \hn(E_s) \,|\,s\in S\}
\subset \ov{C}$, the level set $|S|^{\tau}(E)$ is closed $|S|$.
\end{proposition}

\proof (Sketch) It is enough to assume that $S$ is integral. 
Then by Lemma \ref{one} and noetherian induction, it follows that 
each level subset $|S|^{\tau}(E)$ is constructible. The Lemma \ref{three}.(1)
now implies that the map $s\mapsto \hn(E_s)$ 
is upper semicontinuous. The result follows. \qed

For the following result, see Proposition 3.7 of
Kumar-Narasimhan [Ku-Na]. 

\begin{lemma}\label{five} Let $X$ be a smooth geometrically irreducible 
projective curve over a field $k$, let $E$ be a principal $G$-bundle
on it, and let $(P,\sigma)$ be a reduction of $E$.
Let $N_{\sigma(X), E/P}$ be the normal bundle to the
image $\sigma(X)$ in $E/P$. Then $\sigma^*N_{\sigma(X), E/P}$
is the vector bundle on $X$ associated to the principal $P$-bundle 
$\sigma^*E$ by the representation of $P$ on $Lie(G)/Lie(P)$
induced by the adjoint representation of $G$ on $Lie(G)$. 
If $(P,\sigma)$ is the canonical reduction, then 
$H^0(\sigma(X), N_{\sigma(X), E/P}) = 0$. \hfill$\square$
\end{lemma}

We now come to our main result.

\begin{theorem}\label{six}
{\bf (Main Theorem)} Let $X/S$ be a smooth projective family of 
geometrically irreducible curves over a locally 
noetherian scheme $S/k$, and  
let $E$ be a principal $G$-bundle on $X$. Then for each 
$\tau \in \ov{C}$, there exists a unique locally closed subscheme 
$S^{\tau}(E) \subset S$ with the following universal property.
Any morphism $f: T\to S$ of $k$-schemes factors via 
$S^{\tau}(E) \hra S$ if and only if the pullback $E_T = f^*E$ on $X_T$ 
admits a global relative canonical reduction
$(P, \sigma: X_T \to E_T/P)$ with constant HN type $\tau$ on $T$. 
Moreover, a global relative canonical reduction
for $E_T$, whenever it exists, is unique.  
\end{theorem}

\proof By Proposition \ref{four}, the subset $|S|^{\le \tau}(E)$
is open in $S$. We give it the structure of an open subscheme
of $S$, which we denote as $S^{\le \tau}(E)$. 
The level subset $|S|^{\tau}(E)$ is closed in $S^{\le \tau}(E)$.
In what follows, we
endow $|S|^{\tau}(E)$ with a structure of a closed subscheme $S^{\tau}(E)$
of $S ^{\le \tau}(E)$ which
satisfies the desired universal property. Thus, we can replace 
$S$ by $S^{\le \tau}(E)$ for proving the result, and thereby assume
that $\tau$ is maximal and $|S|^{\tau}(E)$ is closed in $S$.

Let $P = P_{\tau}\subset G$ be the standard parabolic determined by the 
maximal type $\tau$. Let $R_{(E/P)/X/S} \subset Hilb_{(E/P)/S}$ be the
open subscheme given in the statement (B) at the beginning of
Section 2. For each negatively dominant character
$\chi : P\to \GG_m$ of the form $\chi = - \omega_{\alpha}$, 
consider the relatively ample line bundle $L_{\chi}$ 
on $(E/P) \to X$. 
Let $R_{(E/P)/X/S}^{\tau} \subset R_{(E/P)/X/S}$
be the open (and closed) subscheme consisting of all points 
$(P,\sigma)$ (where $\sigma : X_s \to E_s/P$ 
is a reduction) such that $\deg \sigma^*L_{\chi} = \langle \chi, \tau \rangle$. 
By Riemann-Roch, 
this amounts to fixing the Hilbert polynomial w.r.t. each of the finitely many
$L_{\chi}$, so indeed defines an open (and closed) subscheme by 
the requirement of flatness in the definition of $Hilb_{(E/P)/S}$. 

If $R$ is a dvr over $k$ and if 
$\spec R \to S$ is a $k$-morphism under which the generic point maps
to a point in $|S|^{\tau}(E)$, then by closedness of $|S|^{\tau}(E)$
the special point also maps into $|S|^{\tau}(E)$.
The morphism 
$R_{(E/P)/X/S}^{\tau} \to S$ satisfies the valuative 
criterion for properness, as by the Lemma \ref{three} there exist a lift 
$\spec R \to R_{(E/P)/X/S}^{\tau}$. By uniqueness of the canonical
reduction, the morphism $R_{(E/P)/X/S}^{\tau} \to S$ is injective. 
As $R_{(E/P)/X/S}^{\tau}$ is open in $Hilb_{(E/P)/S}$, the
vertical tangent space to $R_{(E/P)/X/S}^{\tau}\to S$ at any point 
$(P,\sigma)$ over $s\in S$ equals the vertical tangent space to 
$Hilb_{(E/P)/S} \to S$ at $[\sigma(X_s)]$. By a standard result on 
relative Hilbert schemes, this equals 
$$H^0(\sigma(X_s), N_{\sigma(X_s), E_s/P})$$ 
where $N_{\sigma(X_s), E_s/P}$
is the normal bundle to $\sigma(X_s) \subset E_s/P$.
By Lemma \ref{five}, this is zero. Also note that the canonical reduction 
$(P,\sigma)$ for any $E_s$ of HN-type $\tau$ exists over $X_s$
itself, that is, $(P,\sigma)$ is $k(s)$-valued (see Remark \ref{can red}). 
The above two facts show
that the projection $R_{(E/P)/X/S}^{\tau} \to S$ is unramified and
induces an isomorphism on residue fields.

Hence we have proved that $\pi: R_{(E/P)/X/S}^{\tau} \to S$ is proper, 
injective, 
unramified, and induces an isomorphism on residue fields $k(\pi(z)) \to k(z)$
for each $z\in R_{(E/P)/X/S}^{\tau}$. Hence by Lemma 4 of [Ni-3],
$\pi$ is a closed embedding. We put $S^{\tau}(E)\subset S$ to be its
image. By its construction, this has the required universal property.
\qed

The Main Theorem has the following immediate corollary.

\begin{corollary}\label{2.7} Under the hypothesis of Theorem \ref{six}, 
if moreover $S$ is reduced and the HN type is globally constant over $S$,
then there exists a relative canonical reduction over $S$.
 \end{corollary}

\bigskip

\newpage

\centerline{\bf The algebraic stacks of HN types}

\medskip

The Main Theorem allows us to show that all principal bundles of a given
HN type form an Artin stack. The proof of this is analogous to
the corresponding result for $\OO$-coherent pure sheaves in [Ni-3],
so we rapidly sketch it.

For any family of curves $X/S$ as above, we have an Artin stack
$Bun(G/X/S)$ over $S$, which to $T\to S$ attaches the groupoid 
$Bun(G/X/S)(T)$ of
all principal $G$-bundles on $X_T$, with pullback (pseudo)-functor
defined the usual way. For any $\tau \in \ov{C}$, let 
$Bun(G/X/S)^{\tau}(T) \subset Bun(G/X/S)(T)$ be the full subgroupoid
of bundles $E_T$ which admit a global relative canonical reduction 
with HN type $\tau$. This condition being preserved under pullbacks
(where you just pull back the reduction), it defines a sub $S$-groupoid
$$Bun(G/X/S)^{\tau} \subset Bun(G/X/S).$$ 
This subgroupoid is actually a stack, as the following effective descent
condition is satisfied.

\begin{lemma}\label{seven} Let $T$ be an $S$-scheme and let 
$E$ be an object of $Bun(G/X/S)^{\tau}(T)$. 
Let $f: T'\to T$ be a faithfully flat quasi-compact morphism. 
If the pullback $f^*E$ is in $Bun(G/X/S)^{\tau}(T')$,
then $E$ is in $Bun(G/X/S)^{\tau}(T)$. 
\end{lemma}

\proof (Sketch) This amounts to saying that the relative 
canonical reduction $(P,\sigma')$ over $T'$ can be descended to $T$. 
Under the two projections 
$\pi_1,\pi_2: T'' \toto T'$, where $T'' = T'\times_TT'$,
the pullbacks of $(P,\sigma')$ are identical 
by the uniqueness of a relative canonical reduction. As a reduction is 
just a section of $E/P$, by the Grothendieck result on effective descent 
for a closed subscheme of a scheme (see Remark \ref{two point nine} below), 
the reduction $(P,\sigma')$ descends to give a reduction $(P,\sigma)$ of $E_T$,
which is canonical as it is pointwise so. \qed

\rem\label{two point nine} (Grothendieck [FGA]). 
Let $X\to T$ be a projective morphism and let 
$f: T'\to T$ be a faithfully flat quasi-compact morphism. 
Let $Y'\subset X_{T'} = X\times_TT'$ be a closed subscheme such that 
under the two projections 
$\pi_1,\pi_2: T'' \toto T'$, where $T'' = T'\times_TT'$,
the two pullbacks of $Y'$ are identical as
closed subschemes of $X_{T''}$. Then there exists a unique
closed subscheme $Y\subset X$ such that $Y' = Y_{T'}$.

\begin{theorem} The stack $Bun(G/X/S)^{\tau}$ is an Artin stack,
which is a locally closed substack of $Bun(G/X/S)$.
\end{theorem}

\proof The Theorem \ref{six} 
implies (just as in the proof of Theorem 8 in [Ni-3])
that the inclusion $1$-morphism $Bun(G/X/S)^{\tau}\hra Bun(G/X/S)$
is representable and it is a locally closed embedding. As 
$Bun(G/X/S)$ is known to be an Artin stack, the result follows. \qed

\section{HN Stratification for families 
of $\Lambda$-modules}

Carlos Simpson introduced in [Si] the notion of $\Lambda$-modules, 
which is a common generalization of important examples such as 
$\OO$-modules, Higgs sheaves, integrable connections, integrable 
logarithmic connections, etc. 
Further, by working over the parameter scheme $S = \AA^1_k$, 
the notion of $\Lambda$-modules incorporates Deligne's idea of continuously 
interpolating between Higgs bundles and integrable connections. 

The result we prove here showing the existence of a schematic
HN stratification for families of $\Lambda$-modules therefore applies to
all the interesting special cases of $\Lambda$-modules that are stated
above, including the interpolation 
between Higgs bundles and integrable connections. They are a generalization of
the results in Nitsure [Ni-3] for $\OO$-modules. The main theorem is given 
together with more expository details in a chapter of the PhD thesis 
(unpublished) of Gurjar [Gu], based on this joint work. 

The subsection 3.1 will recall basic definitions and facts which --
explicitly or essentially --
occur in the literature (mainly in Simpson [Si]) or are `standard knowledge'. 
The main theorem is proved in subsection 3.2.

\subsection{Preliminaries on $\Lambda$-modules}

Let $S$ be a locally noetherian base scheme 
(need not even be equicharacteristic).
Let $f: X\to S$ be a projective, faithfully flat morphism with a chosen
relatively ample line bundle $\OO_{X/S}(1)$. Simpson introduced in 
[Si] the concept of a split almost-polynomial sheaf of rings of 
differential operators $\Lambda$ on $X$ over $S$. 
By definition, such a $\Lambda$ is a sheaf of rings on $X$  
(not necessary commutative), together with a
ring homomorphism $i:\OO_X \rightarrow \Lambda$, a filtration of
$\Lambda$ by subsheaves of abelian groups
$\Lambda_0\subset\Lambda_1\subset \cdots$, and 
a sheaf homomorphism $\zeta : Gr_1(\Lambda) \to \Lambda_0$ 
which satisfies the following properties (1) to (8). 

(1) $\Lambda= \bigcup_{i=0}^{\infty}\Lambda_i$
and $\Lambda_i \cdot \Lambda_j \subset \Lambda_{i+j}$.\\
(2) The image of the given homomorphism $i: \OO_X \rightarrow \Lambda$ is
equal to $\Lambda_0$. (We will denote the resulting homomorphism 
$\OO_X \to \Lambda_0$ again by $i$.)\\
(Note that the above two conditions make each
$\Lambda_i$ a left as well as a right $\OO_X$-submodule of $\Lambda$.)\\
(3) The image of $f^{-1}(\OO_S)$ by the composite 
$f^{-1}\OO_S \stackrel{f^{\sharp}}{\to} \OO_X \stackrel{i}{\to}\Lambda$ 
is contained in the center of $\Lambda$. \\
(4) The left and right $\OO_X$-module structures on 
$Gr_i(\Lambda)= \Lambda_i/\Lambda_{i-1}$ are equal.\\
(5) The sheaves of $\OO_X$-modules $Gr_i(\Lambda)$ are coherent. \\
(6) The sheaf of graded $\OO_X$-algebras $Gr(\Lambda)= \oplus_{i=0}
Gr_i(\Lambda)$ is generated by 
$Gr_1(\Lambda)$ in the sense that the homomorphism
$Gr_1(\Lambda) \otimes_{\OO_X} \cdots \otimes_{\OO_X} Gr_1(\Lambda) 
\to Gr_i(\Lambda)$ is surjective.\\
(7) The condition that $\Lambda$ is {\it almost polynomial}: 
The homomorphism $i: \OO_X \to \Lambda_0$
should be an isomorphism, $Gr_1(\Lambda)$ should be locally free over    
$\OO_{X}$, and the graded ring $Gr(\Lambda)$ should be
naturally isomorphic to the symmetric $\OO_X$-algebra
on $Gr_1(\Lambda)$ under the homomorphism induced by the
multiplication on $\Lambda$.  \\
(8) The condition that $\Lambda$ is {\it split}: 
The given homomorphism $\zeta: Gr_1(\Lambda)\to \Lambda_1$ 
is a homomorphism of left $\OO_X$-modules which splits 
the projection $\Lambda_1 \to Gr_1(\Lambda)$.

By abuse of notation, the data $(\Lambda, \Lambda_j, i, \zeta)$ will
be denoted simply by $\Lambda$.

We now fix an algebra $\Lambda$ as above. 
As Simpson explains in [Si], given any $T \to S$, we functorially get a
split almost-polynomial algebra $\Lambda_T$ on $X_T/T$, which we
call as the {\it pullback of $\Lambda$ under $T\to S$}. Its
underlying sheaf of $\OO_{X_T}$-modules is the pullback
of the left $\OO_X$-module $\Lambda$. The data of the multiplication
on $\Lambda_T$, the filtration $(\Lambda_T)_j$, 
the map $i_T: \OO_{X_T}\to \Lambda_T$
and the splitting $\zeta_T$ are defined functorially.

\medskip

\centerline{\bf $\Lambda$-modules, semistability, HN types}

Let $(Y,\OO_Y(1))$ be a projective scheme over a base field $k$, together
with an ample line bundle $\OO_Y(1)$ on it. Let $\Lambda$ be 
a split almost-polynomial sheaf of rings of 
differential operators on $Y$ over $k$. Unless otherwise indicated,
by a $\Lambda$-module we will
mean a sheaf $E$ of left $\Lambda$-modules on $Y$, which is assumed to be
coherent as an $\OO_Y$-module, 
where the $\OO_Y$-module structure on $E$ is induced via 
$i: \OO_Y \to \Lambda$. Simpson defines a $\Lambda$-modules $E$ on 
$Y$ to be {\it semistable} w.r.t. $\OO_Y(1)$ if 
$E$ is coherent and pure as an $\OO_Y$-module, and moreover for any 
$\Lambda$-submodule $F\subset E$ the inequality 
$r(E)P(F) \le r(F)P(E)$ holds, where for any coherent $\OO_Y$-module 
$M$ whose support is $d$-dimensional, 
we denote by $P(M)(m) = \sum_i (-1)^i\dim_k H^i(Y, M(m))$
the Hilbert polynomial of $M$, and by $r(M)$ the non-negative integer
such that the leading coefficient of $P(M)$ is $r(M)/d!$.  
We put $r(M) =0$ if $M =0$. 

The set $\hnt$ of all abstract HN types for $\OO_Y$-coherent, 
pure $d$-dimensional $\Lambda$-modules on $Y$ is 
the same as the set of all abstract HN types for coherent pure $d$-dimensional 
$\OO_Y$-modules, which we now recall from [Ni-3]. A polynomial 
$f\in \Q[\lambda]$ is called a {\it numerical polynomial} if 
$f(\Z) \subset \Z$. If a nonzero numerical polynomial  
$f$ has degree $d$, it can be uniquely expanded as 
$f = (r(f)/d!)\lambda^d + \mbox{ lower degree terms}$, where
$r(f) \in \Z$. If $f=0$ we put $r(f) =0$.
There is a {\it total order} $\le$ on $\Q[\lambda]$ under which $f\le g$ 
if $f(m) \le g(m)$ for all sufficiently large integers $m$. 
The {\it set of all HN types}, denoted by $\hnt$,
is the set of all 
finite sequences $(f_1,\ldots, f_p)$ of numerical polynomials in 
$\Q[\lambda]$, where $p$ is allowed to vary
over all integers $\ge 1$,  
such that the following three conditions are satisfied.

\noindent{(1)} We have $0< f_1 < \ldots < f_p$ in $\Q[\lambda]$,\\
\noindent{(2)} the polynomials $f_i$ are all of the same degree,
say $d$, and \\
\noindent{(3)} the following inequalities are satisfied  
$$\frac{f_1}{r(f_1)}  > \frac{f_2 - f_1}{r(f_2) - r(f_1)} > \ldots >
\frac{f_p - f_{p-1}}{r(f_p) - r(f_{p-1})}.$$

As described in [Ni-3], to each HN type $(f_1,\ldots, f_p)$
there corresponds a certain subset (HN polygon) 
$$\HNP( f_1,\ldots, f_p) \subset \Z \times \Q[\lambda]$$
which is the union of the segments 
$\ov{x_0x_1} \cup \ov{x_1x_2} \cup \ldots
\cup \ov{x_{p-1}x_p}$
where $x_0 = (0,0)$ and $x_i = (r(f_i), f_i)$ for $1\le i\le p$.
A point $(a,f)\in \Z \times \Q[\lambda]$ is said to {\it lie under}
another point $(b,g )\in \Z \times \Q[\lambda]$ if 
$a=b$ in $\Z$ and $f\le g$ in $\Q[\lambda]$. 
A point $(a,f)\in \Z \times \Q[\lambda]$ is said to {\it lie under the 
polygon} $\HNP(g_1,\ldots, g_q)$ if there exists some $(b,g) \in 
\HNP( g_1,\ldots, g_q)$ such that the point 
$(a,f)$ lies under the point $(b,g)$.
There is a {\it partial order} $\le$ on $\hnt$, under which
$(f_1,\ldots, f_p) \le (g_1,\ldots, g_q)$
if for each $1\le i \le p$, the point 
$(r(f_i), f_i)$ lies under $\HNP(g_1,\ldots, g_q)$.

When a $\Lambda$-module $E$ is $\OO_Y$-coherent of pure dimension $d\ge 0$ 
but not necessarily semistable, 
it admits a unique strictly increasing filtration 
$0= \hn_0(E) \subset\hn_1(E)\subset \ldots \subset \hn_{\ell}(E) = E$
by $\Lambda$-submodules  $\hn_i(E)$ such that for each $1\le i\le \ell$, 
the graded piece $\gr_i(E) = \hn_i(E)/\hn_{i-1}(E)$ is a semistable 
$\Lambda$-module of pure dimension $d$, and the inequalities
$$\frac{P(\gr_1(E))}{r(\gr_1(E))} > \ldots > 
\frac{P(\gr_{\ell}(E))}{r(\gr_{\ell}(E))}$$
hold. This filtration is 
called the {\it Harder-Narasimhan filtration} of $E$
(in the sense of Gieseker semistability). 
The first step $\hn_1(E)$ is called the {\it maximal destabilizing
subsheaf} of $E$.
The integer $\ell$ (also written as $\ell(E)$) is called as the 
{\it length} of the Harder-Narasimhan filtration of $E$. 
In these terms, a nonzero $\OO_Y$-coherent pure-dimensional
$\Lambda$-module $E$ is semistable 
if and only if its Harder-Narasimhan filtration is of length $\ell(E) = 1$.
The ordered $\ell(E)$-tuple 
$$\hn(E) = (P(\hn_1(E)), \ldots ,  P(\hn_{\ell}(E)))\in \hnt$$
is called the {\it Harder-Narasimhan type} of $E$. The reader should note that
the phrases `HN filtration' and 
`HN type' in what follows are to be understood in the sense of 
$\Lambda$-modules.

\rem\label{why the quotient is again pure} 
If $(f_1,\ldots, f_p)\in \hnt$, then $(f_2-f_1, \ldots, f_p-f_1)$
is again in $\hnt$. Let $E$ be an $\OO_Y$-coherent pure-dimensional 
$\Lambda$-module on $Y$ with 
$\hn(E) \le (f_1,\ldots, f_p)\in \hnt$.
If $E'\subset E$ is a $\Lambda$-submodule with $P(E') = f_1$,
then we must have $\hn_1(E) = E'$, that is, such an $E'$
is automatically the maximal destabilizing subsheaf of $E$.
The quotient $E'' = E/E'$ is pure, with 
$\hn(E'') \le (f_2-f_1, \ldots, f_p-f_1)$. Moreover 
$\HOM_{\Lambda}(E', E'') = 0$,  
where $\HOM_{\Lambda}$ denotes the global $\Lambda$-homomorphisms.

\rem\label{HN base changes}
If $(Y,\OO_Y(1))$ is a projective scheme over a field $k$ and if
$K$ any extension field of $k$, then a $\Lambda$-module $E$ on $Y$ is
semistable  w.r.t. $\OO_Y(1)$ if and only if its base-change 
$E_K = E\otimes_kK$ to $Y_K$ is a semistable $\Lambda_K$-module
w.r.t. $\OO_{Y_K}(1) = \OO_Y(1)\otimes_kK$. Consequently, if $E$ is any
$\OO_Y$-coherent pure-dimensional 
$\Lambda$-module on $Y$ then the Harder-Narasimhan filtration
$\hn_i( E_K)$ is just the pullback $\hn_i( E)\otimes_kK$
of the Harder-Narasimhan filtration of $E$.

\medskip

\centerline{\bf Families and quot schemes for $\Lambda$-modules}

Let $f: X\to S$ be a projective, faithfully flat morphism with a chosen
relatively ample line bundle $\OO_{X/S}(1)$, where $S$ is a 
locally noetherian scheme.
Let there be given a split almost-polynomial sheaf of rings of 
differential operators $\Lambda$ on $X$ over $S$.
By a {\it family} of $\Lambda$-modules on $X/S$ 
we will mean an $\OO_X$-coherent
sheaf $E$ of (left) $\Lambda$-modules on $X$, such that as an 
$\OO_X$-module, $E$ is flat over the base $S$. We say that the family
$E$ of $\Lambda$-modules 
is of {\it pure dimension $d$} if moreover each restriction $E|_{X_s}$ 
as an $\OO_{X_s}$-module (where $X_s$ is the schematic fiber over $s\in S$)
is a coherent $\OO_{X_s}$-module of pure dimension $d$. The scheme $S$ is
called the {\it parameter scheme} of the family.

If $E$ is a family of $\Lambda$-modules on $X/S$, and 
$q: E\to F$ is a surjection of $\OO_X$-modules where $F$ is 
a coherent $\OO_X$-module, we say that 
$q$ is a {\it $\Lambda$-quotient} if the kernel of $q$ (which is 
a priori only a coherent $\OO_X$-submodule of $E$) is a $\Lambda$-submodule
of $E$. In that case, note that $F$ acquires the $\Lambda$-module structure
of the quotient $E/\ker(q)$.

\begin{lemma}\label{lambda quotient closed condition} {\rm (Simpson [Si])}
Given any coherent $\OO_X$-module quotient 
$q: E\to F$ of the family $E$ of $\Lambda$-modules on $X/S$, there exists
a unique closed subscheme $S_{\Lambda}\subset S$ which has the 
following universal
property. Any morphism $T\to S$ of schemes factors via 
$S_{\Lambda}$ if and only if
the $\OO_{X_T}$-module pullback $q_T : E_T \to F_T$ on $X_T$ is a $\Lambda_T$
quotient in the above sense.
\end{lemma} 
\proof The proof occurs within the proof of the Theorem 3.8 of [Si]. 
\qed

With $X/S$, $\OO_{X/S}(1)$ and $\Lambda$ as above, let $E$ be a family of 
$\Lambda$-modules on $X/S$ and let $f \in \Q[\lambda]$ be any polynomial.
Let $Q=Quot_{E/X/S}^{f,\OO_{X}(1)}$ be the relative Quot scheme of coherent
$\OO_X$-module quotients of $E_s$ 
with Hilbert polynomial $f$ w.r.t. $\OO_{X/S}(1)$.
By definition, for any $T\to S$, a $T$-valued point of $Q$ is an equivalence
class of quotients as $\OO_{X_T}$-modules $q_T : E_T \to F$ where 
$F$ is a coherent $\OO_{X_T}$-module flat over $T$, and two such 
quotients $q_T : E_T \to F$ and $q'_T : E_T \to F'$ are deemed to be 
equivalent if there exists an $\OO_{X_T}$-module isomorphism $\phi:F\to F'$
such that $q'_T = \phi\circ q_T$. Moreover, it is required that for each
$t\in T$, the Hilbert polynomial of $F_t = F|_{X_t}$ w.r.t. $\OO_{X_t}(1)$ 
is the given polynomial $f \in \Q[\lambda]$, where 
$\OO_{X_t}(1)$ denotes the pullback of $\OO_{X/S}(1)$ to the fiber $X_t$ of
$X_T/T$. The existence of such a scheme $Q$ over $S$, 
and its projectivity over $S$,
are fundamental theorems of Grothendieck (see for example [Ni-1]
for an expository account). As $Q$ represents the above functor,
applying the Yoneda lemma we get a universal family of 
quotients $q: E_Q \to F$ on $X_Q$.

\rem\label{lambda quot scheme}{\bf (Quot scheme for $\Lambda$-modules.)}
With $X/S$, $\OO_{X/S}(1)$ and $\Lambda$ as above, let $E$ be a family of 
$\Lambda$-modules on $X/S$ and let $f \in \Q[\lambda]$ be any polynomial.
Let $Q = Quot_{E/X/S}^{f,\OO_{X}(1)}$ be the relative Quot scheme of coherent
$\OO_X$-module quotients of $E$ with Hilbert polynomial $f$ w.r.t. 
$\OO_{X/S}(1)$.
Let $Q_{\Lambda} \subset Q$ be its closed subscheme defined by Lemma
\ref{lambda quotient closed condition}. By its construction, $Q_{\Lambda}$ 
represents the contravariant functor from $S$-schemes to $Sets$, 
which associates to any $S$-scheme $T$ the set of all $\OO_{X_T}$-coherent
$\Lambda_T$-module quotient of $E_T$ which are flat over $T$ and whose
restriction to each schematic fiber $X_t$ for $t \in T$ has Hilbert 
polynomial $f$ w.r.t. $\OO_{X_t}(1)$. The universal quotient 
$q: E_Q \to F$ on $X_{Q_{\Lambda}}$ will be the restriction (as $\OO$-modules) 
of the universal quotient on $X_Q$.

\begin{proposition}\label{vertical tangent quot scheme}
{\bf (Vertical tangent spaces to the quot scheme $Q_{\Lambda}$)}
Let $k$ be a field, $X$ a projective scheme over $k$,
$\Lambda$ a split almost polynomial sheaf of rings of differential operators 
on $X$, $E$ an $\OO_X$-coherent
$\Lambda$-module, and $q_0 : E \to F_0$ a surjective $\OO_X$-homomorphism
of coherent $\OO_X$-modules such that $q_0$ is a $\Lambda$-quotient as
defined earlier. Let $q_0\in Q_{\Lambda}$ again denote the corresponding
$k$-valued point. Then the tangent space at $q_0$ to $Q_{\Lambda}$ is given by
$$T_{q_0}(Q_{\Lambda}) = \HOM_{\Lambda}(\ker(q_0), F_0).$$
\end{proposition}

\proof Let $Q$ be the quot scheme of coherent $\OO$-quotients
of $E$, and let $Q_{\Lambda}\subset Q$ be its closed subscheme
of $\Lambda$-quotients, as above. 
As proved by Grothendieck, 
the tangent space $T_{q_0}Q$ is given by
$$T_{q_0}Q = \HOM_{\OO_X}(\ker(q_0), F_0)$$
where the right-hand side is the $k$-vector space of all global $\OO_X$-linear
homomorphisms $\ker(q_0) \to F_0$ (see for example [Ni-2] for an exposition). 

Note that any $v\in T_{q_0}Q$ is just a $k[\epsilon]/(\epsilon^2)$-valued 
point of $Q$ which specializes to $q_0$ under $\epsilon \mapsto 0$.
Let $X[\epsilon] = X\otimes_kk[\epsilon]/(\epsilon^2)$, and let
$E[\epsilon]$ be the $\OO$-module pullback of $E$ to $X[\epsilon]$. 
Then in terms of valued points, $v$ can be regarded as a pair 
$(q:E[\epsilon]\to F,\, i: F_0 \to F|_X)$ where $F$ is a coherent 
$\OO_{X[\epsilon]}$-module that is flat over $k[\epsilon]/(\epsilon^2)$,
$q$ is an $\OO_{X[\epsilon]}$-linear surjection, and $i$ is an $\OO_X$-linear
isomorphism. Two such pairs $(q: E[\epsilon] \to F, i)$
and $(q': E[\epsilon] \to F',\, i')$ define the same $v\in T_{q_0}Q$
if and only if there exists
a $\OO_{X[\epsilon]}$-linear isomorphism $F\to F'$ which takes $(q,i)$ to 
$(q',i')$.

As a sheaf of abelian groups 
$E[\epsilon] = E\oplus\epsilon E$ and its germs of local sections
are of the form $(\alpha, \epsilon \beta)$ where 
$\alpha$ and $\beta$ are germs of local sections of $E$. 
Given any $v\in \HOM_{\OO_X}(\ker(q_0), F_0)$, let
$G \subset E[\epsilon]$ be the sub sheaf of abelian groups  
whose germs of local sections are of the form $(\alpha, \epsilon \beta)$
such that $q_0(\alpha) = 0$ and 
$$v(\alpha) = q_0(\beta).$$
Then $G$ is a coherent $\OO_{X[\epsilon]}$-submodule of $E[\epsilon]$
which canonically restricts to $\ker(q_0)$ modulo $\epsilon =0$, and 
the quotient $F = E[\epsilon]/G$ is flat over $k[\epsilon]/(\epsilon^2)$,
with $F_0 = F/\epsilon F$
(flatness of $F$ over $k[\epsilon]/(\epsilon^2)$ is equivalent to the
map $\epsilon: F/\epsilon F \to \epsilon F$ being an isomorphism). Then
the resulting deformation 
$(q: E[\epsilon]\to F, i) \in Q(k[\epsilon]/(\epsilon^2))$ 
represents $v$.

Note that $E[\epsilon]$ is naturally a $\Lambda[\epsilon]$-module.
The point $v\in T_{q_0}Q$ lies in $T_{q_0}(Q_{\Lambda})$ if and only if 
the deformation
$(q: E[\epsilon] \to F, i)$ as concretely described above lies in
$Q_{\Lambda}(k[\epsilon]/(\epsilon^2)) \subset 
Q(k[\epsilon]/(\epsilon^2))$.  
This is clearly the case if and only if $G$ is a 
$\Lambda[\epsilon]$-submodule of $E[\epsilon]$, which  
is equivalent to the condition that $v:  \ker(q_0) \to F_0$
is $\Lambda$-linear. \qed

\bigskip

\centerline{\bf Semicontinuity of HN type}

The well-known Narasimhan-Ramanathan argument 
that semistability is a Zariski open condition for a flat
family of pure dimensional coherent sheaves of $\OO$-modules, which 
uses the existence and projectivity of quot schemes, works equally well
for a pure dimensional family of $\Lambda$-modules when the 
scheme $Q_{\Lambda}$ given by Remark \ref{lambda quot scheme}
is used in place of the quot scheme for $\OO$-quotients. 

Given a family $E$ of $\Lambda$-modules of pure dimension $d$,
we can define a function 
$$|S|\to \HNT: s\mapsto \hn(E_s)$$ 
on the underlying topological 
space $|S|$ of the parameter scheme $S$, taking values in the partially
ordered set $\HNT$ of all HN-types. 
For a family $E$ of $\OO$-modules (which is the special case $\Lambda = \OO_X$),
Shatz proved that the above function is upper semicontinuous. 
Actually, Shatz considered HN types in the sense of 
$\mu$-semistability rather than in the sense of
Gieseker semistability, but his proof readily works for 
HN types in the sense of Gieseker semistability. The key ingredient is again
the existence and projectivity of an appropriate quot scheme for $\OO$-modules.
When we use the schemes $Q_{\Lambda}$ in place of quot schemes, 
the proof of Shatz equally works for our case of 
families of $\Lambda$-modules, with  
HN types in the sense of Gieseker $\Lambda$-semistability 
(this is given in explicit detail in Gurjar [G]). In particular, 
we get the following.

\rem\label{strata are locally closed} 
Let $E$ be a pure dimensional family of $\Lambda$-modules on $X/S$
where $S$ is locally noetherian.
For any $\tau \in \hnt$, the corresponding level set 
$$|S|^{\tau}(E) =\{ s \in |S| \mbox{ such that } \hn(E_s) = \tau \}$$ 
is locally closed in $|S|$, 
the subset 
$|S|^{\le \tau}(E)  = \bigcup_{\alpha \le \tau} |S|^{\alpha}(E) \subset |S|$ 
is open in $|S|$, and $|S|^{\tau}(E)$ is closed in $|S|^{\le \tau}(E)$.

\subsection{Schematic HN stratification}

Let $E$ be a family of pure dimensional $\Lambda$-modules on $X/S$,
with notation as before. A {\it relative HN filtration} for $E$ will mean 
a strictly increasing filtration 
$0= \hn_0(E) \subset\hn_1(E)\subset \ldots \subset \hn_{\ell}(E) = E$
by $\Lambda$-submodules  $\hn_i(E)$ such that for each $1\le i\le \ell$, 
the graded piece $\gr_i(E) = \hn_i(E)/\hn_{i-1}(E)$ 
(regarded as a coherent $\OO_X$-module) is flat over $S$, and
for each $s\in S$, the restriction of the filtration to the fiber $X_s$ 
of $X/S$ (which is a filtration of $E_s = E|_{X_s}$ 
given the flatness condition) is the HN filtration for the $\Lambda_s$-module
$E_s$.

\begin{theorem}\label{main theorem for lambda modules}
Let $X\to S$ be projective, faithfully flat, with a chosen
relatively ample line bundle $\OO_{X/S}(1)$, where $S$ is a locally 
noetherian scheme. Let $\Lambda$ be a split almost-polynomial
sheaf of rings of differential operators on $X/S$.
Let $E$ be a left $\Lambda$-module on $X$, which as an $\OO_X$-module is  
coherent and flat over $S$, and such that the restriction $E_s = E|_{X_s}$
is a pure $\OO_{X_s}$-module
for each $s \in S$. Then we have the following.

{\noindent (1)} Each HN stratum 
$|S|^{\tau}(E)$ of the $\Lambda$-module 
$E$ has a unique structure of a locally
closed subscheme $S^{\tau}(E)$ of $S$, with the following universal
property: a morphism $T\to S$ factors via $S^{\tau}(E)$ if and only if
the pullback $E_T$ on $X\times_ST$ 
admits a relative HN filtration of type $\tau$.

{\noindent (2)} A relative HN
filtration on $E$, if it exists, is unique.

{\noindent (3)} For any morphism $f: T\to S$ of locally noetherian schemes,
the schematic stratum  $T^{\tau}(E_T)\subset T$ for $E_T$ equals the 
schematic inverse image of $S^{\tau}(E)$ under $f$. 
\end{theorem}

\proof With all the preparation that we have made, the proof is now 
exactly the $\Lambda$-module analogue of the proof of the Theorem 5 in 
[Ni-3]. We give a sketch for completeness. 
Let $\tau = (f_1,\ldots, f_{\ell})\in \hnt$ have length $\ell$. 
If $\ell =1$, we take $S^{\tau}(E)$ to be the semi-stable stratum 
$S^{ss}(E)$, which is an open subscheme of $S$. This satisfies all
requirements, as semistability is indeed an open condition on the 
parameter scheme as a special case of the semicontinuity of $\hn(E_s)$, 
and semistability is preserved by 
arbitrary base changes by Remark \ref{HN base changes}. For $\ell \ge 2$,
we proceed by induction. Let $S^{\le \tau}(E) \subset S$ be the open subscheme 
where $\hn(E_s) \le \tau$. Then 
the subset $|S|^{\tau}(E)\subset S^{\le \tau}(E)$ is closed.
We will give it the structure of a closed subscheme 
$S^{\tau}(E)$ of the scheme $S^{\le \tau}(E)$, such that the
conclusion of the theorem is satisfied.
The type $\tau$ is a maximal type for the restriction of the family $E$ to
the inverse image of $S^{\le \tau}(E)$ in $X$.
Thus, it is enough to prove the above theorem for the special case where
$\tau$ is the global maximum type for the given family. 
Hence we now assume that $\tau$ is the global maximum 
type for our given family $E$ on $X/S$. 

Consider the relative quot scheme 
$$Q_{\Lambda}\subset Quot_{E/X/S}^{f_\ell - f_1, \OO_{X/S}(1)}$$ 
of $\OO$-coherent $\Lambda$-quotients of fibers of $E$, 
with Hilbert polynomial $f_\ell - f_1$ where $\tau = (f_1,\ldots, f_{\ell})$.  
Let $q\in Q_{\Lambda}$ be any point, and let $q\mapsto s\in S$.
Let $k(s)\hra k(q)$ be the resulting residue field extension. 
Let $E_q = E|_{X_q}$, and let $q : E_q \to F$ also denote the corresponding 
$\Lambda$-quotient represented by $q$. Then $\ker(q) = 
\hn_1(E_q)$ by Remark \ref{why the quotient is again pure}.
By Remark \ref{HN base changes}
the quotient $q$ is the pullback of the quotient 
$E_s\to E_s/\hn_1(E_s)$ which is defined over $X_s$. 
Hence the residue field extension $k(s) \to k(q)$ is trivial. 
By the uniqueness of $\hn_1(E_s)$, there exists 
at most one such $q$ over $s$.
The fiber of $\pi: Q_{\Lambda} \to S$ over $s$ is the Quot scheme 
$$\pi^{-1}(s) = (Quot_{E_s/X_s/k(s)}^{f_{\ell} - f_1 ,\, \OO_{X_s}(1)})_{\Lambda}$$
of all $\Lambda$-quotients of $E_s$ with Hilbert polynomial $f_{\ell} - f_1$.
By Proposition \ref{vertical tangent quot scheme}, the tangent space
$T_q(\pi^{-1}(s))$ to the fiber $\pi^{-1}(s)$ at $q$ is given by
$$T_q(\pi^{-1}(s)) = \HOM_{\Lambda_s}(\hn_1(E_q),E_q/\hn_1(E_q)).$$ 
This is zero by 
Remark \ref{why the quotient is again pure}.
Hence $\pi: Q_{\Lambda} \to S$ is unramified. 
By Lemma 4 of [Ni-3], any 
morphism $f: T\to S$ between locally noetherian schemes
is a closed embedding if (and only if) 
$f$ is proper, injective, unramified and induces an 
isomorphism $k(f(t)) \to k(t)$ of residue fields for all $t\in T$.
Hence $\pi :  Q_{\Lambda} \to S$ is a closed imbedding.

Now consider the universal $\Lambda$-quotient 
$E_{Q_{\Lambda}} \to E''$ on $X_{Q_{\Lambda}} = X\times_S{Q_{\Lambda}}$. 
By Remark \ref{why the quotient is again pure},
for all $q \in Q_{\Lambda}$ the sheaf $E''_q$ on $X_q$ is
pure-dimensional, with 
$$\hn(E''_q) \le \tau'' = (f_2-f_1,\ldots, f_{\ell} - f_1).$$
In particular, we have $Q_{\Lambda}^{\le \tau''}(E'') = Q_{\Lambda}$.
The Harder-Narasimhan type $\tau''$ has length $\ell -1$, 
hence by induction on the length, the closed subset 
$|Q_{\Lambda}|^{\tau''}(E'')$ of $Q_{\Lambda}$ has the structure of a closed 
subscheme $Q_{\Lambda}^{\tau''}(E'') \subset Q_{\Lambda}$ which has the desired
universal property for $E''$. We regard $Q_{\Lambda}$ as a closed subscheme of 
$S$ via $\pi$, and we finally 
define the closed subscheme $S^{\tau}(E)\subset S$ by
putting  
$$S^{\tau}(E) = Q_{\Lambda}^{\tau''}(E'') \subset Q_{\Lambda} \subset S.$$

By their construction it is clear (as in the proof of Theorem 5 of [Ni-3])
that the subschemes $S^{\tau}(E)$, and the resulting schematic stratification 
of $S$, have the desired properties. \qed

\begin{corollary} 
{\bf (Case of constant HN type over a reduced base)}
Let $X$ be a faithfully flat 
projective scheme over a locally noetherian base scheme
$S$, with a relatively ample line bundle $\OO_{X/S}(1)$
and a split almost polynomial sheaf of rings of differential operators
$\Lambda$ on $X/S$. Let $E$ be a left $\Lambda$-module on $X$, 
which as an $\OO_X$-module is coherent and flat over $S$, and such that 
the restriction $E_s = E|_{X_s}$ is a pure $\OO_{X_s}$-module
for each $s \in S$ of a fixed Harder-Narasimhan type 
$\tau \in \hnt$ as a $\Lambda_s$-module. Suppose moreover that $S$ is 
reduced. Then $S = S^{\tau}(E)$, that is, 
$E$ admits a unique relative Harder-Narasimhan filtration.
\end{corollary}

\centerline{\bf $\Lambda$-modules of HN type $\tau$ form an algebraic stack}

Let $X$ be a faithfully flat 
projective scheme over a locally noetherian base scheme
$S$, with 
a split almost polynomial sheaf of rings of differential operators
$\Lambda$ on $X/S$. For any $S$-scheme $T$, let $\Lambda Coh_{X/S}(T)$
denote the groupoid whose objects are  
all families $E$ of left $\Lambda_T$-module on $X_T$, such that as an 
$\OO_{X_T}$-module $E$ is  
coherent and flat over $T$, and such that the restriction $E_t = E|_{X_t}$
is a pure $\OO_{X_t}$-module
for each $t \in T$. The morphisms in this groupoid are $\Lambda_T$-linear
isomorphisms. For any $S$-morphism $T'\to T$, we have
a natural pullback functor $\Lambda Coh_{X/S}(T)\to \Lambda Coh_{X/S}(T')$,
which makes $\Lambda Coh_{X/S}$ into 
an $S$-groupoid. As in the case of $\OO$-modules
(for which see for example [L-MB] 2.4.4), it can be seen that 
$\Lambda Coh_{X/S}$ is an Artin stack over $S$. 

The Theorem \ref{main theorem for lambda modules}
has the following Corollary, with proofs again as  
in the case of $\OO$-modules in [Ni-3]. We omit the details.

\begin{corollary}
Let $X$ be a faithfully flat 
projective scheme over a locally noetherian base scheme
$S$, with a relatively ample line bundle $\OO_{X/S}(1)$
and a split almost polynomial sheaf of rings of differential operators
$\Lambda$ on $X/S$. Let $\tau$ be any Harder-Narasimhan type. 
Then all families of $\Lambda$-modules on $X/S$ which are flat 
families of coherent pure $\OO$-modules on $X/S$ with a fixed 
Harder-Narasimhan type $\tau$ as $\Lambda$-modules form an 
algebraic stack $\Lambda Coh_{X/S}^{\tau}$ over $S$, which is a 
locally closed substack of the algebraic stack $\Lambda Coh_{X/S}$ of 
all families of $\Lambda$-modules on $X/S$ which are flat 
families of coherent pure $\OO$-modules on $X/S$.
\end{corollary}

It should be noted that whenever a boundedness theorem holds for 
semistable $\Lambda$-modules, each stack $\Lambda Coh_{X/S}^{\tau}$
admits an atlas of finite type over $S$, as in Proposition 9 of [Ni-3].

\bigskip

\bigskip

\parskip=4pt

{\large \bf References}

[A-B] Atiyah, M.F. and Bott, R.: The Yang-Mills equations over 
Riemann surfaces. 
Philos. Trans. Roy. Soc. London Ser. A 308 (1983), no. 1505, 523–-615.

[Be] Behrend, K.: Semi-stability of reductive group schemes over curves. 
Math. Ann. 301 (1995), no. 2, 281–-305.

[Bi-Ho] Biswas, I. and Holla, Y.: 
Harder-Narasimhan reduction of a principal bundle. 
Nagoya Math. J. 174 (2004), 201–-223.

[Gu] Gurjar, S. : {\it Topics in principal bundles}. 
PhD thesis, Tata Institute of Fundamental Research, 2012.

[L-MB] Laumon, G. and Moret-Bailly, L. : {\it `Champs alg\'ebriques.'}
Springer (2000).

[Ku-Na] Kumar, S. and Narasimhan, M.S. : Picard group of the moduli
spaces of $G$-bundles. Math. Ann. 308 (1997), no. 1, 155-173.

[Ni-1] Nitsure, N. : Construction of Hilbert and Quot schemes. Part
2 of {\it `Fundamental Algebraic Geometry -- Grothendieck's FGA
Explained.'},  Fantechi et al, Math. Surveys and Monographs Vol.
123, American Math. Soc. (2005).

[Ni-2] Nitsure, N. : Deformation theory for vector bundles. Chapter 5
of {\it Moduli Spaces and Vector Bundles} 
(edited by Brambila-Paz, Bradlow, Garcia-Prada and Ramanan), 
London Math. Soc. Lect. Notes 359,  
Cambridge Univ. Press (2009).

[Ni-3] Nitsure, N. : 
Schematic Harder-Narasimhan stratification. 
Internat. J. Math. 22 (2011), no. 10, 1365–1373. 

[Sh] Shatz, S.S. : The decomposition and specialization of
algebraic families of vector bundles. Compositio Math. 35
(1977), no. 2, 163--187.

[Si] Simpson, C. Moduli of representations of the fundamental 
group of a smooth projective variety -I.  
~~Publ. Math. I.H.E.S. 79  (1994), 47--129.

\bigskip

\bigskip

Sudarshan Gurjar \hfill e-mail: {\tt sgurjar@math.tifr.res.in}\\
Nitin Nitsure  \hfill e-mail: {\tt nitsure@math.tifr.res.in}\\
School of Mathematics, Tata Institute of Fundamental Research,\\
Homi Bhabha Road, Mumbai 400 005, India.\\

\centerline{27-Aug-2012}

\end{document}